\documentclass[12pt]{amsart}
\newif\ifdraft\draftfalse %%% This will produce line numbers on everything
\usepackage{amssymb}
\ifdraft\usepackage{lineno}\else\fi

\usepackage{amsrefs}  %%%% Use this one last to avoid complains about things 
\input xypic
\xyoption{all}

\define\lam{\Lambda^{-1,-1}}
\renewcommand{\l}{\lambda}

\define\R{{\Bbb R}}

\define\Z{{\Bbb Z}}
\define\C{{\Bbb C}}
\define\half{\frac{1}{2}}

\renewcommand{\gg}{\mathfrak g}
\define\ad{\text{ad}\,}
\define\Ad{\text{Ad}\,}

\theoremstyle{plain}
\newtheorem{thm}[equation]{Theorem}

\newtheorem{conjecture}[equation]{Conjecture}

\theoremstyle{definition}

\theoremstyle{remark}
\newtheorem{rk}[equation]{Remark}

\numberwithin{equation}{section}

\def\CC{\mathbb C}

\def\GG{\mathbb G}

\def\RR{\mathbb R}

\def\ZZ{\mathbb Z}
\def\Z{\ZZ}
\def\R{\RR}

\def\C{\CC}

\def\ff{\mathcal F}
\def\gg{\frak g}
\def\hh{\mathcal H}

\def\mm{\mathcal M}

\def\uu{\mathcal U}
\def\vv{\mathcal V}
\def\ww{\mathcal W}

\def\zz{\mathcal Z}

\def\Cal{\mathcal}

\def\blnm{\ifdraft\begin{linenomath*}\else\fi}
\def\elnm{\ifdraft\end{linenomath*}\else\fi}

\def\M{\mm}
\def\V{\vv}

\def\Ydag{Y^{\ddag}}

\def\SL{\mathop{\mathrm{SL}}}
\def\GL{\mathop{\mathrm{GL}}}

\def\Ext{\mathop{\mathrm{Ext}}}
\def\AV{\mathop{\mathrm{AV}}\nolimits}

\def\ad{\mathop{\mathrm{ad}}\nolimits}
\def\Ad{\mathop{\mathrm{Ad}}\nolimits}
\def\Gr{\mathrm{Gr}}

\begin{document}

\title[Zero locus]{The zero locus of an admissible normal function}
\author{Patrick Brosnan}
\address{Department of Mathematics\\
The University of British Columbia\\
Room 121, 1984 Mathematics Road\\
Vancouver, B.C., Canada V6T 1Z2}
\email{brosnan@math.ubc.ca}
\author{Gregory J. Pearlstein}
\address{Mathematics Department
Duke University, Box 90320\\
Durham, NC 27708-0320}
\email{gpearlst@math.duke.edu}

\subjclass[2000]{32G20, 14D07,14D05}
\date{}

\def\lineme{\ifdraft\linenumbers\else\fi}
\lineme
\begin{abstract}
We prove that the zero locus of an admissible normal
function over an algebraic parameter space $S$ is algebraic in the
case where $S$ is a curve.  
\end{abstract}

\maketitle
\section{Introduction}

Let $S$ be a smooth, complex projective variety.  Following Morihiko
Saito~\cite{MS}, we define an admissible normal function on $S$ to be
an admissible variation of graded-polarized mixed Hodge structure $\vv$
on the complement of a normal crossing divisor $D\subseteq S$ which is
an extension of the trivial variation $\ZZ(0)$ by a variation of pure,
polarized Hodge structure $\hh$ of weight $w<0$.  That is an admissible
normal function is an element $\nu\in\Ext^1_{\AV}(\ZZ(0),\hh)$ where
$\AV$ denotes the abelian category of admissible variations of mixed Hodge
structure on $S$ which are smooth on $S-D$.

Henceforth, we assume that $w=-1$.  In this case, an admissible normal
function corresponds to the usual notion of a horizontal normal
function on $S$ with moderate growth along $D$ together with existence
of a suitable relative weight filtration along each irreducible
component of $D$.  In this article (Theorem~\ref{thm-main}), we settle
the following conjecture communicated to us by M.~Green and
P.~Griffiths in the case where $S$ is a curve.

\begin{conjecture}
  Let $\nu$ be an admissible normal function on $S$.  Then the zero
  locus $\zz$ of $\nu$ is an algebraic subvariety of $S$.
\end{conjecture}

In analogy with~\cite{CDK}, an unconditional proof of this conjecture
provides indirect evidence in support of standard conjectures on
higher regulators and filtrations on Chow groups.

A rough outline of our proof is a follows: Let $\uu$ be an open subset
of $S$ in the analytic topology which does not intersect $D$.  Then
the zero locus of $\nu$ on $\uu$ is complex analytic since the
restriction of $\nu$ to $\uu$ is a holomorphic section of associated
bundle of intermediate Jacobians.  Thus, in order to prove that the
zero locus of $\nu$ is algebraic, it is sufficient to show that:
\begin{itemize}
\item[(*)] For each point $p\in D$ there exists an analytic open neighborhood 
$\uu_p \subset S$ of $p$ on which the zero 
locus of  $\nu$ has only finitely many components.
\end{itemize}
In the case where $S$ is a curve, we  verify $(*)$ using the orbit
theorems of the second author and results of P.~Deligne.  

The canonical real grading $Y(s)$ (described below) of the Hodge
structure $\vv_s$ at a point $s\in S-D$ will play a crucial role in
our proof.  The central idea is that $\nu$ is $0$ at $s$ if and only
if $Y(s)$ is integral.  It is therefore crucial to understand the
asymptotics of $Y(s)$ as $s$ tends to a point $s_0\in D$.  In
Theorem~\ref{thm:3.5}, we use the $\SL_2$-orbit theorem of~\cite{P3}
to show that $\Ydag:=\lim_{s\to s_0} Y(s)$ exists for $s_0\in D$.
Now, it is clear that $\nu$ can only vanish in a neighborhood of $s_0$
if $\Ydag$ is non-integral.  Knowing that the limit exists allows us
to concentrate on the case where $\Ydag$ is not integral.  This case
can then be handle by a rather explicit computation of the zero locus
in the neighborhood of $s_0$.

\subsection*{Acknowledgments} We gratefully acknowledge helpful conversations
with P.~Deligne, M.~Green, P.~Griffiths, J.~Lewis and S.~Usui.  We would also
like to thank the Institute for Advanced Study in Princeton, New Jersey
for hosting us (in nearby offices) while most of this work was done.

\section{The Zero locus at a smooth point}
%%%
%%% Start Here
%%%

\par As a preliminary step in our analysis of the zero locus of $\nu$
at infinity, we derive the local defining equations
of $\mathcal Z$ at an interior point of $S$.  To this end, we begin
with a review of mixed Hodge structures and their gradings.

\subsection*{Gradings} Let $V$ be a finite dimensional vector space
over a field $k$.  A grading for $V$ is a direct sum
decomposition $V=\oplus V_k$ of $V$ into subspace $V_k$ indexed by
integers.  An $n$-grading of $V$ is a direct sum decomposition
$V=\oplus V_{w}$ indexed by $n$-tuples of integers.  It is well-known
(and easy to see) that gradings are in one-one correspondence with
linear actions of the multiplicative group $\GG_m$ on $V$, and
$n$-gradings are in one-one correspondence with linear actions of the
$n$-torus $\GG_m^n$ on $V$: To an $n$-grading $V=\oplus V_w$ one
associates the action where $(t_1,\ldots, t_n)v= t_1^{w_1}\cdots
t_n^{w_n}v$.  Conversely, to an action of $\GG_m^n$ on $V$ one obtains
a grading by writing $V$ as a direct sum of its isotypical subspaces.

If $k$ is a field of characteristic $0$, then  gradings are in
one-one correspondence with semi-simple endomorphisms $Y$ of $V$ with
integral eigenvalues.  The correspondence is the one that takes a $\GG_m$ 
action
on $V$ to its derivative at $1\in\GG_m$ viewed as an endomorphism of $V$.   
Conversely, it takes an endomorphism $Y$ to the direct sum decomposition 
$V=\oplus_{k\in\ZZ} V_{Y,k}$ where $V_{Y,k}=\{v\in V|Yv=kv\}$.
Similarly, $n$-gradings are in one-one correspondence with $n$-tuples 
of commuting semi-simple endomorphisms $Y_1,\ldots, Y_n$ with integral
eigenvalues.

Throughout this paper, the vector spaces $V$ which will occur will be
over fields of characteristic $0$.  We will work with the first notion
and the last notion of a grading interchangeably.

Suppose that $V$ is equipped with a filtration 
\blnm
\begin{equation}
  V=W_r V\supset W_{r-1} V \supset \cdots \supset W_{l-1}V=\{0\}
\end{equation}
\elnm
by subspaces $W_i$ with $i\in\ZZ$.  A {\em grading} of $W_{\bullet}$ is then
a grading $V=\oplus_{i\in\ZZ} V_i$ of $V$ such that
$W_k=\oplus_{i\leq k} V_i$.

\subsection*{Deligne's Grading}  
We now recall a fundamental result of Deligne concerning mixed Hodge 
structures. (See \cite{Hodge2}*{Lemme 1.2.8} for the original theorem or 
\cite{CKS}*{Theorem 2.13}
where the result appears in the notation used below.)

A mixed Hodge structure $(F,W)$ induces a functorial bigrading 
\begin{equation}
\label{eq:2.1}
         V_{\C} = \bigoplus_{p,q}\, I^{p,q}		
\end{equation}
such that
\begin{enumerate}
\item $F^p = \bigoplus_{a\geq p}\, I^{a,b}$;
\item $W_k = \bigoplus_{a+b\leq k}\, I^{a,b}$;
\item $\bar I^{p,q} = I^{q,p} \mod \oplus_{r<q,s<p}\, I^{r,s}$.
\end{enumerate}
In particular, a mixed Hodge structure $(F,W)$ induces a grading of $W$
via the semisimple endomorphism $Y:V_{\C}\to V_{\C}$ which acts as 
multiplication by $p+q$ on $I^{p,q}$.  We will call this grading 
{\em Deligne's grading} $Y_{F,W}$.

\subsection*{Normal functions}
Returning now to the normal function setting, let $S$ be a complex
manifold of dimension $n$ and $\vv$ be the variation of mixed Hodge 
structure associated to $\nu$. Let $p\in\ZZ$ and $(s_1,\dots,s_n)$ be 
local holomorphic coordinates on a polydisk $\Delta^n\subseteq S$ which
vanish at $p$.  Then, since $\Delta^n$ is simply connected, we can parallel  
translate the data of $\vv$ back the reference fiber $V = \vv_p$.  The Hodge 
filtration $\ff$ of $\vv$ then corresponds to a holomorphic, horizontal 
decreasing filtration $F(s)$ of 
$V_{\CC}$. The weight filtration $\ww$ of $\vv$ corresponds 
to a constant filtration $W$ of $V_{\ZZ}$ with weight graded-quotients
\blnm
\begin{equation*}
          \Gr^W_0(V_{\ZZ}) = \ZZ(0),\qquad \Gr^W_{-1}(V_{\ZZ}) = H_{\ZZ}
\end{equation*}
\elnm
and $\Gr^W_k = 0$ for $k\neq 0,-1$.  Similarly, the graded-polarizations
of $\ww$ correspond to constant polarizations of $\Gr^W$.  

On account of the short length of $W$, the grading
\blnm
\begin{equation}
\label{eq:2.2}
          Y(s) = Y_{(F(s),W)}					
\end{equation}
\elnm
defined by \eqref{eq:2.1} can be characterized as the unique real grading of 
$W$ which preserves $F(s)$.  If $Y_{\ZZ}$ is any integral grading of $W$ then 
the image $1\in\Z(0)$ under the map
\blnm
$$
        Y(s) - Y_{\Z} :\ZZ(0) \to H_{\RR}/H_{\ZZ}
$$
\elnm
gives the point in the Griffiths intermediate Jacobian corresponding to
the extension \eqref{eq:2.1} via the standard isomorphism
\blnm
$$
        H_{\R}/H_{\Z} \cong \frac{H_{\C}}{F^0 H_{\C} + H_{\Z}}.
$$
\elnm
Accordingly, $p$ belongs to the zero locus of $\nu$ if and only if
$Y(p)$ is an integral grading of $W$.  Consequently, since $Y(s)$
is real analytic\footnote{The grading $Y(s)$ is holomorphic when
viewed as a section of the bundle of intermediate Jacobians}
in $s$ and the set of integral gradings of $W$ is a
discrete subset of the affine space of $\R$-gradings of $W$, there 
exists a neighborhood of $p$ in which the zero locus of $\nu$ is given 
by the equation
\blnm
$$       
       Y(s) = Y(p).                                                    
$$
\elnm

 The filtration $F(s)$ takes its values in a classifying space $\M$ 
of graded-polarized mixed Hodge  structure~\cite{P1}.  Let $G_{\C}$ denote the
Lie group consisting of all automorphisms of $V_{\C}$ which preserve 
$W$ and act by complex isometries on $\Gr^W$.  Then, for each point $F\in\M$
there exists a neighborhood $U_{\C}$ of zero in the Lie algebra $\gg_{\C}$
such that the map
\blnm
\begin{equation}
\label{eq:2.3}
        u \mapsto e^u.F                                         
\end{equation}
\elnm
is a holomorphic submersion from $U_{\C}$ onto a neighborhood of $F$
in $\M$.  Accordingly, if 
\begin{equation}
  \gg_{\C} = \bigoplus_{p+q\leq 0}\, \gg^{p,q}       \label{eq:lie-bigrading}
\end{equation}
denotes the Deligne bigrading of the induced mixed Hodge structure 
$(F\gg_{\C},W\gg_{\C})$ then map $(2.3)$ restricts to a biholomorphism 
from a neighborhood of zero in subalgebra
\blnm
$$
        q_F = \bigoplus_{p<0,p+q\leq 0}\, \gg^{p,q}
$$   
\elnm 
onto a neighborhood of $F$ in $\M$.

 Letting $F=F(p)$, the constructions of the previous paragraph show
that near $p$ we can write 
\blnm
$$
     F(s) = e^{\Gamma(s)}.F
$$
\elnm
relative to a unique holomorphic function $\Gamma(s)$ with values in 
$q_F$ which vanishes at $p$.  Let $Y=Y(p)$ and 
\blnm
$$
     \Gamma(s) = \Gamma_0(s) + \Gamma_{-1}(s)
$$
\elnm
denote the decomposition of $\Gamma(s)$ according to the eigenvalues of 
$Y$.  By the Campbell--Baker--Hausdorff formula, there exists a universal
power series $\Psi(t)$ such that
\blnm $$
     e^{u+v}e^{-u} = e^{\Psi(\ad u)v}.
$$ \elnm
In particular, 
\blnm $$
\aligned
     e^{\Gamma(s)}.Y 
        &= e^{\Gamma_0(s) + \Gamma_{-1}(s)}e^{-\Gamma_0(s)}e^{\Gamma_0(s)}.Y \\
        &= e^{\Psi(\ad\Gamma_0(s))\Gamma_{-1}(s)}.Y                     
         = Y + \Psi(\ad\Gamma_0(s))\Gamma_{-1}(s)
\endaligned
$$ \elnm
is a holomorphic grading of the weight filtration (over $\C$) which preserves 
$F(s)$.  Therefore, there exists a real analytic section $\zeta(s)$ of 
$\gg_{\C}^{F(s)}\cap W_{-1}\gg_{\C}$ such that
\blnm $$
       Y(s) = Y + \Psi(\ad\Gamma_0(s))\Gamma_{-1}(s) + \zeta(s)
$$ \elnm
and hence the equation $Y(s) = Y(p)$ is equivalent to
\begin{equation}
\label{eq:2.4}
       \Psi(\ad\Gamma_0(s))\Gamma_{-1}(s) + \zeta(s)  = 0.      
\end{equation}
Equation \eqref{eq:2.4} implies that, near $p$ on the zero locus of $\nu$,
\begin{equation}
\label{eq:2.5}
        \Psi(\ad\Gamma_0(s))\Gamma_{-1}(s)\in 
               \gg_{\C}^{F(s)}\cap W_{-1}\gg_{\C}              
\end{equation}
on the zero locus of $\nu$.  Conversely, whenever equation \eqref{eq:2.5}
 holds,
$Y=Y(p)$ is a real grading of $W$ which preserves $F(s)$.  Because these
two properties specify $Y(s)$ uniquely, it then follows that whenever
equation \eqref{eq:2.5} holds, $Y(s) = Y(p)$.  Thus, on a neighborhood of $p$,
the zero locus of $\nu$ is given by equation \eqref{eq:2.5}.

 Applying $e^{-\ad\Gamma(s)}$ to both sides of \eqref{eq:2.5}, it then follows 
that the equation for the zero locus near $p$ is 
\begin{equation}
\label{eq:2.6}
       e^{-\ad\Gamma(s)}\Psi(\ad\Gamma_0(s))\Gamma_{-1}(s) 
                \in \gg_{\C}^F\cap W_{-1}\gg_{\C}.                 
\end{equation}
However, $q_F$ is a nilpotent subalgebra of $\gg_{\C}$ which is a vector
space complement to $\gg_{\C}^F$ in $\gg_{\C}$.  Furthermore, $q_F$ is
closed under $\ad Y$.  Therefore, $\Gamma(s)$, $\Gamma_0(s)$ and 
$\Gamma_{-1}(s)$ take values in $q_F$, and
$
        e^{-\ad\Gamma(s)}\Psi(\ad\Gamma_0(s))\Gamma_{-1}(s) 
$
takes values in $q_F$.  Consequently, equation \eqref{eq:2.6} is equivalent to
\begin{equation}
\label{eq:2.7}
        e^{-\ad\Gamma(s)}\Psi(\ad\Gamma_0(s))\Gamma_{-1}(s) = 0.
\end{equation}
Thus, since $e^{-\ad\Gamma(s)}$ is invertible, equation \eqref{eq:2.7} implies
the following result.

\begin{thm}
\label{thm:2.8} 
Near $p$, the zero locus of $\nu$ is given by the equation
$\Gamma_{-1}(s)=0$.
\end{thm}
\begin{proof} Applying $e^{\ad\Gamma(s)}$ to \eqref{eq:2.7} implies
that the zero locus is given by the equation
\begin{equation}
\label{eq:2.9}
        \Psi(\ad\Gamma_0(s))\Gamma_{-1}(s) = 0                
\end{equation}
By the Campbell--Baker--Hausdorff formula,
$$ 
    \Psi(u)v = v + \sum_{j>0}\, c_j (\ad u)^j v
$$
and hence
\begin{equation}
        \Psi(\ad\Gamma_0)\Gamma_{-1} 
        = \Gamma_{-1} + \sum_{j>0}\, c_j (\ad\Gamma_0)^j \Gamma_{-1}.
        \label{eq:2.10}
\end{equation}
Consequently, if 
$$
         \Gamma_0 = \sum_{k>0}\, \Gamma^{-k,k},\qquad         
         \Gamma_{-1} = \sum_{\ell>0}\, \Gamma^{-\ell,\ell-1}
$$
denote the decomposition of $\Gamma_0$ and $\Gamma_{-1}$ into Hodge
components with respect to the bigrading \eqref{eq:lie-bigrading}
then
$$
      \Psi(\ad\Gamma_0)\Gamma_{-1} 
      = \Gamma^{-1,0} \mod \bigoplus_{r\geq 2}\, \gg^{-r,r-1}.
$$
As such, the equation \eqref{eq:2.9} then implies that $\Gamma^{-1,0} = 0$.  
Proceeding by induction, assume that $\Gamma^{-\ell,1-\ell} = 0$ for $\ell<n$.
Then,
$$
      \Psi(\ad\Gamma_0)\Gamma_{-1} 
      = \Gamma^{-n,n-1} \mod \bigoplus_{r\geq n+1}\, \gg^{-r,r-1}
$$
and hence the equation \eqref{eq:2.9}, $\Gamma^{-n,n-1} = 0$.  Thus, 
$\Gamma_{-1}=0$ is the local defining equation for $\mathcal Z$.
\end{proof}

\section{Limiting Grading} In this section, we prove that
when $S$ is a curve, the grading $(2.2)$ has a well defined limit
$Y^{\ddag}$ as $s$ approaches a puncture $p\in S$.  Simple examples
show that in higher dimensions, the limiting value of $(2.2)$ depends
not only on the point in the boundary divisor but also the direction
of approach.  

Let $\Delta\subset S$ be a disk containing the puncture $p$.
By passing to a finite cover if necessary, we can assume that 
the local monodromy of the restriction of $\V$ to the punctured
disk $\Delta^* = \Delta - \{p\}$ is unipotent.  Let $s$ be a 
local coordinate on $\Delta$ which vanishes at $p$, let $A$
be an angular sector of $\Delta^*$ and $s_o$ be a point in $A$.
Then, we can parallel translate the Hodge filtration of $\V$ back
to a single valued filtration $F(s)$ on $V=\V_{s_o}$.  Analytic 
continuation of $F(s)$ to all of $\Delta^*$ then gives the multivalued 
period map
\blnm $$
        \varphi:\Delta^*\to\Gamma\backslash\M
$$ \elnm
of $\V$.  By local liftablity, there exists a holomorphic, horizontal lifting
of $\varphi$ to a map $\tilde F$ from the upper half-plane $U$ into $\M$
which makes the following diagram commute.
\blnm $$
\xymatrix{
U\ar[r]^{\tilde F}\ar[d]_{e^{2\pi iz}} &  \M\ar[d]           \\
         \Delta^*\ar[r]^{\varphi} & \Gamma\backslash\M     \\
}
$$ \elnm

Furthermore, upon picking a branch of $\log(s)$ on $A$ and letting 
\blnm $$
        z = x + iy = \frac{1}{2\pi i}\log(s)
$$ \elnm
there is unique lifting $\tilde F(z)$ such that, for $s\in A$,
\blnm $$
        \tilde F(z) = F(s).
$$ \elnm

By unipotent monodromy, we $\tilde F(z+1) = e^N.\tilde F(z)$ and hence
\blnm $$
        \tilde\varphi(z) = e^{-zN}.\tilde F(z)
$$ \elnm
drops to a map $\tilde\varphi$ from $\Delta^*$ into the 
\lq\lq compact dual" of $\M$.  The admissibility of $\V$ then asserts that 
\begin{itemize}
\item[{(a)}] $F_{\infty} = \lim_{s\to 0}\tilde\varphi(s)$ exists;
\item[{(b)}] The relative weight filtration $M$ of $W$ and $N$ exists.
\end{itemize}
From these properties, together with Schmid's nilpotent orbit theorem,
Deligne then deduces \cite{SZ} that the pair $(F_{\infty},M)$ is a 
mixed Hodge structure relative to which $N$ is a $(-1,-1)$-morphism.

The mixed Hodge structure $(F_{\infty},M)$ induces a mixed Hodge
structure on $\gg_{\C}$ with Deligne bigrading 
\begin{equation}
\label{eq:3.1}
        \gg_{\C} = \bigoplus_{r,s}\gg^{r,s}.                     
\end{equation}
Note that in equation $(3.1)$, it is possible to have $r+s>0$ since elements
of $\gg_{\C}$ only preserve $W$ and not $M$.  Nonetheless, the nilpotent
subalgebra
\blnm $$
        q_{\infty} = \bigoplus_{r<0}\,\gg^{r,s}                  
$$ \elnm
is a vector space complement to $\gg_{\C}^{F_{\infty}}$.  Reasoning as
in \S 2 (cf. [P1]), it then follows that near the puncture 
$s=0$ we can write
$
        \tilde\varphi(s) = e^{\Gamma(s)}.F_{\infty}
$
relative to a unique holomorphic function $\Gamma(s)$ which takes values
in $q_{\infty}$ and vanishes at $s=0$.  Untwisting the definition of 
$\tilde\varphi$, it then follows that 
\begin{equation}
\label{eq:3.2}
    F(s) = e^{\frac{1}{2\pi i}\log(s) N}e^{\Gamma(s)}.F_{\infty}       
\end{equation}
over the angular sector $A$.

To determine the asymptotic behavior of the grading 
\blnm $$
        Y(s) = Y_{(F(s),W)}
$$ \elnm
on $A$ we shall use equation \eqref{eq:3.2} together with the $SL_2$-orbit theorem of 
\cite{P3} and a result of Deligne which constructs a grading $Y$
of the weight filtration $W$ which is well adapted to both $N$ and
the limiting mixed Hodge structure $(F_{\infty},M)$.  

More precisely, suppose that $Y_M$ is a grading of $M$ which 
preserves $W$ and satisfies
\blnm $$
        [Y_M,N] = -2N.
$$ \elnm
Then, Deligne (see \cite{KP}*{Appendix, Theorem 1})
shows that there exists a unique, functorial grading
\begin{equation}
\label{eq:3.3}
        Y' = Y'(N,Y_M)                                  
\end{equation}
such that $Y'$ commutes with both $N$ and $Y_M$\footnote{The general
statement~\cite{P2} of Deligne's result for longer weight filtrations 
involves the interplay between the decomposition of $N$ according to $\ad Y'$
and the graded representations of $sl_2$.}. Furthermore,
\begin{enumerate}
\item[(a)] If $Y_M$ is defined over $\R$ the so is $Y'$;
\item[(b)] If $(F,M)$ is a mixed Hodge for which $N$ is a $(-1,-1)$-morphism
           and induces sub mixed Hodge structures on $W$ then the grading
           $Y'$ produced from $N$ and the grading of $M$ by the $I^{p,q}$'s
           of $(F,M)$ preserves $F$. 
\end{enumerate}

To show the existence of the limiting grading 
\blnm $$
        Y^{\ddag} = \lim_{s\to 0}\, Y(s)
$$ \elnm
we now invoke the $SL_2$-orbit theorem of \cite{CKS}*{Proposition 2.20}:  Let 
\blnm $$
        (\hat F,M) = (e^{i\delta}.F_{\infty},M)
$$ \elnm
denote Deligne's splitting of the limiting mixed Hodge structure of $\V$
and 
\blnm $$
        \lam = \bigoplus_{r,s<0}\, \gg^{r,s}_{(\hat F,M)}.
$$ \elnm
Define $G_{\R} = G_{\C}\cap \GL(V_{\R})$ and let $\gg_{\R}$ denote the
Lie algebra of $G_{\R}$.  Then, there exists a distinguished, real 
analytic function
\blnm $$
        g:(a,\infty):(a,\infty)\to G_{\R}
$$ \elnm
and element
\blnm $$      
        \zeta \in \gg_{\R}\cap\ker(N)\cap\lam
$$ \elnm
such that
\begin{enumerate}
\item[(a)] $e^{iyN}.F_{\infty} = g(y)e^{iyN}.\hat F$;
\item[(b)] $g(y)$ and $g^{-1}(y)$ have convergent series expansions about
$\infty$ of the form
\blnm $$
\aligned
       g(y) &= e^{\zeta}(1 + g_1 y^{-1} + g_2 y^{-2} + \cdots)      \\
  g^{-1}(y) &= (1 + f_1 y^{-1} + f_2 y^{-2} + \cdots)e^{-\zeta}         
\endaligned
$$ \elnm
with $g_k$, $f_k\in \ker(\ad N)^{k+1}$;
\item[(c)] $\delta$, $\zeta$ and the coefficients $g_k$ are related by
the formula
\blnm $$
      e^{i\delta} 
      = e^{\zeta}\left(1 + \sum_{k>0}\, \frac{1}{k!}(-i)^k(\ad\, N_0)^k\,g_k
                 \right).
$$ \elnm
\end{enumerate}

Combining this result with equation \eqref{eq:3.2}, we obtain the following
asymptotic formula for $Y(s)$ over the angular sector $A$: 
\blnm $$
\aligned
        F(s)  = e^{zN}e^{\Gamma(s)}.F_{\infty}   
             &= e^{xN}e^{\Gamma_1(s)}e^{iyN}.F_{\infty}         \\
             &= e^{xN}e^{\Gamma_1(s)}g(y)e^{iyN}.\hat F     
              = e^{xN}g(y)e^{\Gamma_2(s)}e^{iyN}.\hat F
\endaligned
$$ \elnm
where $\Gamma_1(s) = \Ad(e^{iyN})\Gamma(s)$ and                          
$\Gamma_2(s) = \Ad(g^{-1}(y))\Gamma_1(s)$.

 Let $\hat Y_M$ denote the grading of $M$ defined by the $I^{p,q}$'s of 
$(\hat F,M)$ and $\hat Y$ be the grading of $W$ defined by application
of Deligne's construction to the pair $(N,\hat Y_M)$.  Then \cite{P3},
\blnm $$
        H = \hat Y_M - \hat Y
$$ \elnm
belongs to $\gg_{\R}$ and satisfies $[H,N] = -2N$.  Furthermore, since
$\hat Y_M$ and $\hat Y$ preserve $\hat F$, so does $H$.  Therefore,
\blnm $$
        e^{iyN}.\hat F = y^{-\half H}.F_o
$$ \elnm
where $F_o = e^{iN}.\hat F$.   By the $SL_2$-orbit theorem, $F_o$
belongs to $\M$.  Consequently,
\blnm $$
        F(s) = e^{xN}g(y)e^{\Gamma_2(s)}y^{-\half H}.F_o
             = e^{xN}g(y)y^{-\half H} e^{\Gamma_3(s)}.F_o       
$$ \elnm
where 
\blnm $$
        \Gamma_3(s) = \Ad(y^{\half H})\Gamma_2(s)
                    = \Ad(y^{\half H} g(y) e^{iyN})\Gamma(s).      
$$ \elnm

To continue, observe that since,
$
        y = -\frac{1}{2\pi}\log|s|
$
and $H$ has only finitely many eigenvalues (all of which are integral), the 
action of  $\Ad(y^{\half H})$ on $\gg_{\C}$ is bounded by an integral
power of $y^{\half}$. Similarly, since $g(y)$ is bounded as $s\to 0$, so
is the action of $\Ad(g(y))$.  Likewise, since $N$ is nilpotent, the
action of $\Ad(e^{iyN})$ on $\gg_{\C}$ is bounded by a power of $y$.  
Therefore, since $\Gamma(s)$ is a holomorphic function of $s$ which
vanishes at $s=0$, $\Gamma_3(s)$ is a real analytic function on $A$
which satisfies the growth condition
\blnm $$
        \Gamma_3(s) = O((\log|s|)^b s)
$$ \elnm
for some half integral power $b$.  In particular, near $s=0$,
\blnm $$
        Y_{(e^{\Gamma_3(s)}.F_o,W)} = Y_{(F_o,W)} + \gamma_4(s)
$$ \elnm
for some real analytic function $\gamma_4(s)$ which is again of
order $(\log|s|)^b s$.  By Deligne \cite{CKS},
\blnm $$
        Y_{(F_o,W)} = \hat Y.
$$ \elnm
Therefore,
\blnm $$
\aligned
        Y(s) &=  e^{xN}g(y)y^{-\half H}.Y_{(e^{\Gamma_3(s)}.F_o,W)}   \\
             &=  e^{xN}g(y)y^{-\half H}.(Y+ \gamma_4(s))             \\     
             &=  e^{xN}g(y).(\hat Y + \gamma_5(s))
\endaligned
$$ \elnm
where $\gamma_5(s) = \Ad(y^{-\half H})\gamma_4(s)$ is again
of order $\log|s|^{b'} s$ for some half-integral power $b'$.  

 Define
\blnm $$
\aligned
        \tilde g(s) 
        &= e^{xN}g(y)e^{-xN}                        \\
        &= e^{\zeta}\left(1 + \sum_{k>0} (\Ad(e^{xN}g_k))y^{-k}\right).
\endaligned
$$ \elnm
Then, since $x= \frac{1}{2\pi}\text{Arg}(s)$ is bounded on the angular
sector $A$,
\blnm $$
        \lim_{s\to 0} g(s) = e^{\zeta}.
$$ \elnm
Consequently, because $N$ commutes with $\hat Y$,
\begin{equation}
\label{eq:3.4}
        Y(s) = \tilde g(s).(\hat Y + Ad(e^{xN})\gamma_5(s)).     
\end{equation}
Therefore, since $\gamma_5(s)$ is order $(\log(s))^{b'} s$, we
can take the limit of equation \eqref{eq:3.4} to obtain:

\begin{thm}
\label{thm:3.5}
\begin{equation}
\label{eq:3.6}
        Y^{\ddag} = \lim_{s\to 0}\, Y(s) = e^{\zeta}.\hat Y.  
\end{equation}
\end{thm}

\begin{rk} Since the right hand side of $(3.6)$ depends 
only on the triple $(F_{\infty},W,N)$, $Y^{\ddag}$ is independent of choice of angular
sector $A$.  Likewise, a change of local coordinate $s$ changes 
$F_{\infty}$ to $e^{\l N}.F_{\infty}$.  Therefore, due to the 
functorial nature of Deligne's construction of the grading $Y'$
and the fact that $[Y',N] = 0$, the right hand side of $(3.6)$ is
independent of the choice of coordinate $s$.  Likewise, since
the right hand side of $(3.6)$ commutes with $N$, it is well defined
independent of the choice of reference fiber.  Consequently, in the
geometric setting, $Y^{\ddag}$ should have a direct geometric meaning.
\end{rk}

\section{Zero Locus at Infinity} To verify the conjecture in the 
case where $S$ is a curve, we now note that the finiteness condition $(*)$ 
is preserved under passage to finite covers.  Therefore, we may assume as
in \S 3 that the associated variation of mixed Hodge structure $\V$ has 
unipotent monodromy about each point $p\in D$.  The requirement that the 
zero locus of $\nu$ has only finitely many irreducible components on a 
neighborhood of $p\in D$ is then equivalent to the existence of a disk 
$\Delta\subset S$ such that $\Delta\cap D = \{p\}$ on which the zero locus of 
$\nu$ is either
\begin{enumerate}
\item[(a)] The empty set;
\item[(b)] All of $\Delta$, in which case $\V$ is the trivial extension
of $\Z(0)$ by $\Cal H$.
\end{enumerate}

 Applying Deligne's construction \eqref{eq:3.3} to the limiting mixed Hodge 
structure $(F_{\infty},M)$, we get a grading $Y_{\infty}$ of $W$ which 
preserves $F_{\infty}$.  Therefore,
\blnm $$
        Y_{\infty}(s) = e^{\frac{1}{2\pi i}\log(s)N}e^{\Gamma(s)}.Y_{\infty}
$$ \elnm
is a (complex) grading of $W$ which preserves the Hodge filtration of 
$F(s)$ near $s=0$ over the angular sector $A$.  Decomposing $\Gamma(s)$ as 
\blnm $$
        \Gamma(s) = \Gamma_0(s) + \Gamma_{-1}(s)
$$ \elnm
according to the eigenvalues of $\ad Y_{\infty}$, it then follows that
\begin{align}
        Y_{\infty}(s) 
        &= e^{\frac{1}{2\pi i}\log(s) N}
                e^{\Psi(\ad\Gamma_0(s))\Gamma_{-1}(s)}.Y_{\infty}\notag\\
        &= e^{\frac{1}{2\pi i}\log(s) N}.
                (Y_{\infty} + \Psi(\ad\Gamma_0(s))\Gamma_{-1}(s))\notag\\
        &= Y_{\infty} + e^{\frac{1}{2\pi i}\log(s) \ad N} 
                        \Psi(\ad\Gamma_0(s))\Gamma_{-1}(s).\label{eq:4.1}
\end{align}

As in the derivation of equation \eqref{eq:2.9}, we then have
\begin{equation}
\label{eq:4.2}
        Y(s) = Y_{\infty}(s) + \zeta(s)                          
\end{equation}
for some section $\zeta(s)$ of $W_{-1}\gg_{\C}\cap\gg_{\C}^{F(s)}$.

 Now, unlike a normal function over $\Delta^n$ considered earlier, the
the function $\zeta(s)$ defined by equation \eqref{eq:4.2} may in principle have
singularities at $s=0$.  To show that this is not the case, observe
that since $\Gamma(s)$ is holomorphic and vanishes at $s=0$ and $N$
is nilpotent,
\begin{equation}
\label{eq:4.3}
        \lim_{s\to 0} e^{\frac{1}{2\pi i}\log(s) \ad N} 
                        \Psi(\ad\Gamma_0(s))\Gamma_{-1}(s) = 0.   
\end{equation}
Therefore, since the limit 
\blnm $$
         Y^{\ddag} = \lim_{s\to 0} Y(s)                
$$ \elnm
exists by Theorem $(3.5)$, equations \eqref{eq:4.2} and \eqref{eq:4.3}
imply that $\zeta(s)$ also has a continuous extension to $0$ in the
angular sector $A$.

In particular, if $Y^{\ddag}$ is not an integral grading of $W$
then there is a neighborhood of zero in angular sector $A$ on which
$Y(s)$ is not integral, and hence $\nu$ has no zeros on this neighborhood.
Thus, it remains to consider the case where $Y^{\ddag}$ is integral.
By \cite{P3}, 
\blnm $$
        \hat Y = e^{-i\delta}.Y_{\infty}
$$ \elnm
and hence
\blnm $$
        Y^{\ddag} =  e^{\zeta}.\hat Y = e^{\zeta}e^{-i\delta}.Y_{\infty}.       
$$ \elnm
We can write 
\blnm $$
        e^{\zeta}e^{-i\delta} = e^{\xi} 
$$ \elnm
for some (unique)
\blnm $$
        \xi \in \ker(\ad N)\cap\lam_{(\hat F,M)}
$$ \elnm
since both $\zeta$ and $\delta$ belong to the subalgebra
$\ker(N)\cap\lam_{(\hat F,M)}$.

 To continue, note that
\blnm $$
       \gg^{r,s}_{(F_{\infty},W)}   
        = e^{i\ad \delta}(\gg^{r,s}_{(\hat F,M)})
$$ \elnm
and hence
\blnm $$
        \lam_{(F_{\infty},W)} 
        =  e^{i\ad \delta}\lam_{(\hat F,M)}
        = \lam_{(\hat F,M)}
$$ \elnm
since $\lam_{(\hat F,M)}$ is closed under $\ad\delta$.  As such
\blnm $$
        \xi\in\ker(\ad N)\lam_{(\hat F,M)} 
        = \ker(\ad N)\cap\lam_{(F_{\infty},M)}.
$$ \elnm
Consequently,
\blnm $$
        Y^{\ddag} = e^{\xi}.Y_{\infty}
                  = Y_{\infty} + \Psi(\ad\xi_0)\xi_{-1}
$$ \elnm
where $\xi_0$ and $\xi_{-1}$ both belong to 
$\ker(\ad N)\cap\lam_{(F_{\infty},M)}$ since $Y_{\infty}$ commutes with $N$
and preserves each summand $\gg^{r,s}_{(F_{\infty},M)}$ of 
$\lam_{(F_{\infty},M)}$.  As such,
\blnm $$
        \Psi(\ad\xi_0)\xi_{-1}\in\ker(\ad N)\cap\lam_{(F_{\infty},M)}.
$$ \elnm
        
 Returning now to equation \eqref{eq:4.2}, it then follows that
\blnm $$
        Y(s) = Y^{\ddag} - \Psi(\ad\xi_0)\xi_{-1} 
                + e^{\frac{1}{2\pi i}\log(s)\ad N}
                  \Psi(\ad\Gamma_0(s))\Gamma_{-1}(s) + \zeta(s) 
$$ \elnm
where $\zeta(s)$ is a real analytic section of 
$\gg_{\C}^{F(s)}\cap W_{-1}\gg_{\C}$.  In particular, since
\blnm $$
        \lim_{s\to 0}\, Y(s) = Y^{\ddag}
$$ \elnm
is integral, it then follows from the continuity of $Y(s)$ that near
$s=0$ the zeros of $\nu$ occur where
\blnm $$
       -\Psi(\ad\xi_0)\xi_{-1} 
                + e^{\frac{1}{2\pi i}\log(s)\ad N}
                  \Psi(\ad\Gamma_0(s))\Gamma_{-1}(s) + \zeta(s) = 0.
$$ \elnm
Equivalently,
\begin{align}
        \Ad (e^{\frac{1}{2\pi i}\log(s) N}e^{\Gamma(s)})^{-1}
        &\left(\Psi(\ad\xi_0)\xi_{-1} 
        - e^{\frac{1}{2\pi i}\log(s)\ad N}\Psi(\ad\Gamma_0(s))\Gamma_{-1}(s)
                  \right)\notag \\
        &= \Ad (e^{\frac{1}{2\pi i}\log(s) N}e^{\Gamma(s)})^{-1} \zeta(s).
           \label{eq:4.4}
\end{align}                                                  

Thus, since the right hand side of equation \eqref{eq:4.4} takes values in 
$\gg_{\C}^{F_{\infty}}$ and the left hand side of \eqref{eq:4.4} takes values
in $q_{\infty}$, it then follows that the zeros of $\nu$ occur exactly
where
\blnm $$
        e^{\frac{1}{2\pi i}\log(s)\ad N}\Psi(\ad\Gamma_0(s))\Gamma_{-1}(s)
        = \Psi(\ad\xi_0)\xi_{-1}.
$$ \elnm
Since $\Psi(\ad\xi_0)\xi_{-1}\in\ker(\ad N)$, the equation can be further 
reduced to just
\blnm $$
      \Psi(\ad\Gamma_0(s))\Gamma_{-1}(s)
        = -\Psi(\ad\xi_0)\xi_{-1}.
$$ \elnm
Because $\Gamma(s)$ is a holomorphic function which vanishes at zero, the
above equation only has solutions near $s=0$ only if 
\blnm $$
        \Psi(\ad\xi_0)\xi_{-1} = 0
$$ \elnm
(i.e. $Y^{\ddag} = Y_{\infty}$).  In this case, the equation is just
\blnm $$
        \Psi(\ad\Gamma_0(s))\Gamma_{-1}(s) = 0.
$$ \elnm
Again, because $\Gamma(s)$ is holomorphic at $s=0$, the solutions to the
above equation are either isolated or all of $A$. 

Thus, we have obtained the following.

\begin{thm}
\label{thm-main}
  Let $\nu$ be an admissible normal function on a complex, projective
curve $S$ smooth outside of a finite set $D\subset S$.  Then the zero locus
$\zz$ of $\nu$ is an algebraic subset of $S-D$.
\end{thm}

\end{document}